\documentclass[12pt]{amsart} \usepackage{amsmath}
\usepackage[margin=1.2in]{geometry}
\usepackage{booktabs} \usepackage{pdfpages}
\usepackage{fancyhdr}

\begin{document}

\newcommand{\Kappa}{{\mathcal{K}}}
\newcommand{\calR}{{\mathcal{R}}}
\newcommand{\Ham}{H}
\newcommand{\dhab}{\Ham_a}
\newcommand{\dhbb}{\Ham_b}
\newcommand{\dhpb}{\Ham_p}
\newcommand{\dhqb}{\Ham_q}
\newcommand{\dhub}{\Ham_u}
\newcommand{\dhvb}{\Ham_v}
\newcommand{\dhxb}{\Ham_x}
\newcommand{\dhyb}{\Ham_y}

\newcommand{\mutwothree}{\mu_{2,3}} 
\newcommand{\muonethree}{\mu_{1,3}}
\newcommand{\muonetwo}{\mu_{1,2}}
\newcommand{\dee}{\delta}

\newcommand{\grad}{\nabla}
\newcommand{\calM}{{\mathcal M}}
\newcommand{\calO}{{\mathcal O}}
\newcommand{\calo}{{\mathcal O}}
\newcommand{\calE}{{\mathcal E}}
\newcommand{\calTR}{{\mathcal TR}}
\newcommand{\PP}{{\mathbb {P}}}
\newcommand{\CC}{{\mathbb {C}}}
\newcommand{\RR}{{\mathbb {R}}}
\newcommand{\TT}{{\mathbb {T}}}
\newcommand{\quot}{\mathop{\rm quot}\nolimits}
\newcommand{\quott}{\mathop{\rm quot} \mathop{\rm quot}\nolimits}
\newcommand{\diag}{\mathop{\rm diag}\nolimits}
\newcommand{\CM}{\mathcal M}

\newcommand{\gradH}{\nabla \Ham}

\newcommand{\ga}{T_a}
\newcommand{\gp}{T_p}
\newcommand{\gq}{T_q}
\newcommand{\gu}{T_u}
\newcommand{\gv}{T_v}
\newcommand{\gx}{T_x}
\newcommand{\gy}{T_y}
\newcommand{\gb}{T_b}

\newcommand{\hab}{d_{\ga}\Ham}
\newcommand{\hbb}{d_{\gb}\Ham}
\newcommand{\hpb}{d_{\gp}\Ham}
\newcommand{\hqb}{d_{\gq}\Ham}
\newcommand{\hub}{d_{\gu}\Ham}
\newcommand{\hvb}{d_{\gv}\Ham}
\newcommand{\hxb}{d_{\gx}\Ham}
\newcommand{\hyb}{d_{\gy}\Ham}

\newcommand{\npb}{\|\gp\|^2}
\newcommand{\nqb}{\|\gq\|^2}
\newcommand{\nub}{\|\gu\|^2}
\newcommand{\nvb}{\|\gv\|^2}
\newcommand{\nxb}{\|\gx\|^2}
\newcommand{\nyb}{\|\gy\|^2}

\newcommand{\ca}{c_{a}}
\newcommand{\cb}{c_{b}}
\title[Spectral curves for the triple reduced product]{Spectral curves for the triple reduced product of coadjoint orbits for $\mbox{SU}(3)$}
\author[J. Hurtubise]{Jacques Hurtubise}\thanks{JH was partially supported by a grant from NSERC}
\address{Department of Mathematics \& Statistics \\
McGill University \\ Montr\'eal, Qu\'ebec \\ Canada}
\email{hurtubis@math.mcgill.ca}
\author[L. Jeffrey]{Lisa Jeffrey}\thanks{LJ was partially supported by a grant from NSERC}
\address{Department of Mathematics \\
University of Toronto \\ Toronto, Ontario \\ Canada}
\email{jeffrey@math.toronto.edu}
\author[S. Rayan]{Steven Rayan}\thanks{SR was partially supported by a New Faculty Grant from the University of Saskatchewan}
\address{Department of Mathematics \& Statistics\\
University of Saskatchewan \\ Saskatoon, Saskatchewan \\ Canada}
\email{rayan@math.usask.ca}
\author[P. Selick]{Paul Selick}\thanks{PS was partially supported by a grant from NSERC}
\address{Department of Mathematics \\
University of Toronto \\ Toronto, Ontario \\ Canada}
\email{selick@math.toronto.edu}
\author[J. Weitsman]{Jonathan Weitsman}\thanks{JW was partially supported by NSF grant DMS-12/11819}
\address{Department of Mathematics \\
Northeastern University \\ Boston, MA 02115\\ USA}
\email{j.weitsman@neu.edu}

\maketitle

\begin{center}May 5, 2017\end{center}

\vspace{5pt}

\begin{center}\emph{Dedicated to Nigel Hitchin on the occasion of his 70th birthday.}\end{center}

\begin{abstract}We give an identification of the triple reduced product of three coadjoint orbits in $\mbox{SU}(3)$ with a space of Hitchin pairs over a genus $0$ curve with three punctures, where the residues of the Higgs field at the punctures are constrained to lie in fixed coadjoint orbits.  Using spectral curves for the corresponding Hitchin system, we identify the moment map for a Hamiltonian circle action on the reduced product.  Finally, we make use of results of Adams, Harnad, and Hurtubise to find Darboux coordinates and a differential equation for the Hamiltonian.\end{abstract}

\section{Introduction}

We consider the symplectic quotient of a product of three coadjoint orbits. 
Let $\lambda$, $\mu$, $\nu$ be
 diagonal $3  \times 3$ traceless matrices with real eigenvalues, so that
 $i\lambda, i\mu, i\nu \in {\mathfrak{su}(3)}$.  Then, let $\calo_{i\lambda}$,
 $\calo_{i\mu}$, $\calo_{i\nu}$ be their corresponding orbits under the
 adjoint action on $\mathfrak{su}(3)$.  We define the \emph{triple reduced product} to
 be the quotient $$\PP(\lambda, \mu, \nu) := (\calo_{i\lambda }\times
 \calo_{i\mu} \times \calo_{i\nu})//\mbox{SU}(3).$$ Here, $\mbox{SU(3)}$ acts diagonally on
 the product of orbits, via the adjoint action.  The notation $//$
 indicates that we are taking the symplectic quotient $M^{-1}(0)/\mbox{SU}(3)$ of
 the product of the three orbits, where $$M(X,Y,Z)=X+Y+Z$$ is the moment
 map for the diagonal adjoint action of the group.  This quotient can be interpreted as a moduli space of triangles embedded in $\mathfrak{su}(3)^*$, with vertices lying in fixed coadjoint orbits. 
  In this way, $\PP(\lambda, \mu, \nu)$ is reminiscent of polygon spaces \cite{KM}.

 \newcommand{\MM}{\mathcal{M}}
 \newcommand{\LL}{\mathcal{L}}

We are interested in $\PP(\lambda, \mu, \nu)$ because the orbit method --- that is, the study of the orbits of the coadjoint action on the dual of
a  Lie algebra --- has wide-ranging applications in geometry. This technique originates in the work of Kirillov (e.g. \cite{KV}).
The symplectic quotient of products of orbits by the diagonal action 
is a prototype for moduli spaces of flat connections on punctured $2$-spheres, where the
holonomy around each puncture is restricted to a specific conjugacy class. Passing to 
the language of holomorphic bundles, these are moduli spaces of parabolic bundles on a genus $0$ surface with $n$ punctures.  When the weights at each puncture are sufficiently small, the corresponding moduli space of parabolic
bundles is identified with the $n$-fold reduced product --- see 
Theorem 6.6 in  \cite{LJ}, for example.
The case $n=3$, which is the first nontrivial case, has a symplectic volume that is given by the work of Suzuki and Takakura
\cite{ST}.  Their formula for the volume is piecewise linear in $\lambda,\mu, \nu$.

Topologically, $\PP(\lambda, \mu, \nu)$ is a $2$-sphere as shown recently in \cite{TRP}.  The main objective in \cite{TRP} was
to find a Hamiltonian function on
 $\PP(\lambda, \mu, \nu)$ whose Hamiltonian flow generates an
 $S^1$-action on it.  We were able to construct such a Hamiltonian function, 
using an auxiliary function $f: \PP(\lambda,\mu,\nu) \to [0,1]$  on $\PP(\lambda,\mu, \nu)$ 
which is  surjective. It is required that 
 level sets of the boundary points $0$  and $1$ are points and the level sets of interior points
are circles. 
 This means that taking the union of the level sets of $f$  identifies
$\PP(\lambda, \mu, \nu)$  with $S^2$.
 However, we could
only define the Hamiltonian function indirectly  as a definite integral
involving $f$.

Of course, there could be many such functions; the question remains of trying to obtain such a Hamiltonian in a fairly natural way. This is where our work on $\PP(\lambda, \mu, \nu)$ makes contact with the world of Hitchin systems.  In the present paper, we show that a Hamiltonian generating an $S^1$ action 
on $\PP(\lambda, \mu, \nu)$ arises naturally from the Hitchin map of a Hitchin system; more properly, the generalized Hitchin systems, developed in \cite{AHH, Ma, Bo, BiRa} and others.

 The core of our approach in this paper is to symplectically identify the $\PP(\lambda, \mu, \nu)$ with a compact moduli space of Higgs
 bundles over an appropriate punctured base, with residues constrained to fixed coadjoint orbits.  The philosophy here is that, in
 trying to study a tuple of matrices, it is often easier to combine
 them into a single object, namely a ``Higgs field'', which will be a
 matrix-valued polynomial with coefficients coming from the fixed orbits.  This was used to great effect in, for instance, the additive Deligne-Simpson problem (as developed in \cite{CRA,DEL,SPN,SBN} and other related works).   Given the identification of certain moduli spaces of parabolic bundles with $n$-fold reduced products as in \cite{LJ}, and given that polygon spaces appear within tame parabolic Hitchin systems \cite{GM,FR}, it is not entirely surprising that one can view $\PP(\lambda, \mu, \nu)$ as a Hitchin-like system.  The identification here, however, is not with a subvariety of a parabolic Hitchin system, as the Higgs bundles corresponding to triangles in $\PP(\lambda, \mu, \nu)$ do not have fixed parabolic structures at the punctures.  Rather, we fix the orbits that the residues lie in.

The natural curve for the associated Hitchin system is a triple-punctured projective line.   We embed $\PP(\lambda, \mu, \nu)$ in a moduli space of $\mbox{SU}(3)$ Higgs bundles on this curve, and this embedding lies in the fixed point set of an antiholomorphic involution. In this way, $\PP(\lambda, \mu, \nu)$ is identified with a compact, real Hitchin system.  It turns out that the desired Hamiltonian function is the Hitchin map --- the map sending a Higgs field to its characteristic polynomial --- for this system.   The Hamiltonian $S^1$ action to which the Hitchin map corresponds rotates these circles at a common speed.

This algebro-geometric expression of our finite-dimensional coadjoint orbits can also be viewed as an embedding of coadjoint orbits into a loop algebra, and indeed this version predates Hitchin's work; see e.g. Adler and van Moerbeke \cite{AvM,AHH}, Reyman and Semenov-Tian-Shansky \cite {RST}, and Mischenko and Fomenko \cite{MF}; the  paper \cite{AHH3} explains the link 
of Adler-Kostant-Symes flows with the flows of Mischenko and Fomenko \cite{MF}.  See also Donagi and Markman \cite{DM},
which explains the link between Adler-Kostant-Symes and generalized Hitchin systems. Either way, the geometric expression one obtains has certain advantages: for example, one sees that the fixed points of our action correspond to singular spectral curves. It also raises questions about what can be done in a more general case.  In higher rank one has, quite easily, an integrable system, but there is a natural question of how this can be turned into an action of a torus --- in other words, of obtaining action-angle variables.

It may be profitable to exploit a possible relationship with quiver varieties which occurs in the case where the coadjoint orbits of compact groups are replaced with the coadjoint orbits of the corresponding complex group.  We recall that closures of complex coadjoint orbits of $\mbox{GL}(n,\CC)$ are Nakajima quiver varieties for a certain kind of framed Dynkin graph \cite{KP0,CRA2,Na}. Building from this, the symplectic quotient of a product of these orbits by the diagonal action of $\mbox{GL}(n,\CC)$ is a so-called ``star-shaped'' quiver variety \cite{CRA}.   In this framework, the manifold supporting our integrable system is a real locus in the additive part of the $\mbox{E}6$ example in, for instance, Section 5 of \cite{Boa} and again in \cite{Boa2,Boa3}.  One might then investigate flows on real loci in the corresponding $\mbox{E}7$ and $\mbox{E}8$ examples.   (We thank Philip Boalch for bringing this idea to our attention, pointing out these references, and explaining the potential connection to $\mbox{E}$-type graphs.)

In any case, we hope that the working out of this simple example echoes, however faintly, some of the recurrent features of Nigel Hitchin's work, with theory working itself out in some quite concrete and pretty algebraic geometry, in particular involving elliptic functions and elliptic curves.

\section{Identification with an $\mbox{SU}(3)$ Hitchin system}

\subsection{Spectral curves and the Hitchin system}

Now, let $\CM$ stand for the moduli space of isomorphism classes of $\mbox{SU}(3)$ Higgs bundles on $\PP^1$ punctured at $z=0,\pm1$.  This means that each point in $\CM$ is a Higgs bundle of the form $(P,\Phi)$ where $P$ is a holomorphic principal $\mbox{SU}(3)$-bundle and$$\Phi\in\mbox H^0(\PP^1,\mbox{ad}(P)\otimes K(D)),$$ where $K\cong\mathcal O(-2)$ is the canonical line bundle on $\PP^1$ and $D=0+(-1)+1$ is the divisor of the marked points, considered as a formal sum.  As a result of taking values in $\mbox{ad}(P)$, the matrix $\Phi(z)$ is trace-free and anti-hermitian for each $z\in\PP^1$.

Next, we restrict to the locus $\CM_0\subset\CM$ where $P$ is isomorphic to the trivial bundle $\PP^1\times\mbox{SU}(3)$.  Along this locus, we need only keep track of $\Phi$, which takes the form$$\Phi(z)=\left(\frac{X}{z}+\frac{Y}{z-1}+\frac{Z}{z+1}\right)dz,$$ where the residues $X$, $Y$, and $Z$ are matrices in $\mathfrak{su}(3)$ and $X+Y+Z=0$.  Morally, there is a stability condition that one must impose when defining $\CM$ --- this is a version of Hitchin's slope stability condition adapted to the structure group and to the divisor $D$ --- but since we restrict to Higgs fields for the trivial bundle, every $\Phi$ is semistable. Furthermore, owing to the gauge-theoretic origins of Higgs bundles, we had defined $\Phi$ as a twisted one-form, but this will be unnecessary for our purposes and so we simply consider each Higgs field to be of the form$$\Phi(z)=\frac{X}{z}+\frac{Y}{z-1}+\frac{Z}{z+1}.$$

Finally, we add the following constraint: $X$, $Y$, and $Z$ must lie along the coadjoint orbits $\calo_{i\lambda}$, $\calo_{i\mu}$, and $\calo_{i\nu}$, respectively.  Such a constraint is required if one wants a symplectic, as opposed to a Poisson, structure.  We call this space $\CM_0(\lambda,\mu,\nu)$. This means that $\CM_0(\lambda,\mu,\nu)$ is the symplectic quotient $M^{-1}(0)/\mbox{SU}(3)$, with $M$ defined as in Section 1, and where $\mbox{SU}(3)$ acts on $\Phi$ by conjugation.  It is clear that $$\CM_0(\lambda,\mu,\nu)\cong\PP(\lambda, \mu, \nu)$$ as symplectic manifolds.  From now on, we simply write $\PP(\lambda, \mu, \nu)$ for both.

There is a map from $\PP(\lambda, \mu, \nu)$ to an affine space $\mathbb A$, called the \emph{Hitchin map}, that takes $\Phi(z)$ to its characteristic polynomial.  As it is more convenient to clear the poles, we instead take the characteristic polynomial of $L(z) = z(z^2-1)\Phi(z)$.  Accordingly, we set 
$$\rho(z,\eta) = \det (z(z^2-1)\Phi(z) -\eta I) = \det (L(z)-\eta I)$$
and one has 
$$L(z) =  (z^2-1) X + z(z+1) Y + z(z-1) Z = (Y-Z) z + (Y+Z)$$
 after taking the moment map condition $X+Y+Z=0$ into account.  The \emph{spectral curve} $S$ is then defined by $\rho(z,\eta)=0$. A quick adjunction calculation shows that the genus of the spectral curve is $1$.
 
 This spectral curve is fixed above $z= 1,0,-1$, as we are keeping $X,Y,Z$ in fixed coadjoint orbits. In fact, setting 
 $$ \rho(z,\eta) = (iH z(z^2-1) + Q_0(z)) + Q_1(z) \eta   - \eta^3,$$
 with $Q_i$ quadratic  functions (there is no $\eta^2$ term, as the matrices have trace zero), one has that the coefficient of the $Q_i$ are fixed (are Casimirs) in $\PP(\lambda, \mu, \nu)$, and in fact the only coefficient of the Hitchin map which is nonconstant along the orbit is $iH$. This function, up to additive constants along $\PP(\lambda, \mu, \nu)$, can be taken to be $$-i\det(Y-Z)$$ or $$-i{\rm {tr}}(Y-Z)^3/3.$$

\subsection{Lax form for $H$}
A function $H$ gives a Hamiltonian flow along $\PP(\lambda, \mu, \nu)$, which can be written in Lax form (See \cite{Hi, AHH}). If one identifies $dH$ with an element of the polynomial loop algebra on ${\mathfrak{su}(3)}$, to be thought of as acting on elements $f$ of the tangent space of $\PP(\lambda, \mu, \nu)$ consisting of ${\mathfrak{su}(3)}$-valued functions with poles at $-1,0,1$ by the trace residue pairing $<dH,f >= {\rm Tr}\,\mbox{Res}_\infty(dH f)$, the flow is given by 
$$ \dot \Phi = [ dH, \Phi].$$
 In our case, $H = i{\rm {tr}}(Y-Z)^3/3$ gives us
 $$dH = i(z^{-1} L(z)^2)_+ = i((Y-Z)^2 z + (Y-Z)(Y+Z) + (Y+Z) (Y-Z)).$$
 Accordingly, our flow becomes 
 $$L'(z) = i[ (Y-Z)^2 z + (Y-Z)(Y+Z) + (Y+Z) (Y-Z), (Y-Z) z +  (Y+Z)],$$
 giving 
 $$ (Y-Z)' = 0$$
 $$ (Y+Z)' = i [ (Y+Z) (Y-Z),   (Y+Z)].$$ 
 Hence, we have
 $$ Y' = i[Y, YZ+ZY + Z^2]$$
 and a similar equation for $Z'$, showing that $Y,Z, Y+Z$ all evolve by conjugation, and so one stays inside the moduli space $\PP(\lambda, \mu, \nu)$.
 
Indeed, the matrix $L(z)$ as a whole is evolving by conjugation (by a $z$-dependent matrix), and the spectral curve is thus a constant of motion.  What is flowing is an extra geometric datum: a sheaf supported on the spectral curve. Indeed,  over the total space  of the line bundle $\calo(1)  \to \PP^1$  we have
an exact sequence
$$0\rightarrow \calo^{\oplus 3}(-1) \buildrel{L(z)-\eta I}\over{\longrightarrow} \calo^{\oplus 3} \longrightarrow \calE$$
 that defines a sheaf $ \calE$ supported over the spectral curve. When the curve is smooth and reduced, one can show that  $ \calE$ is a line bundle on the curve. The sheaf $\calE$, along with the curve, encodes $\Phi$.  Indeed, the fiber of the Hitchin map on the full moduli space over the complex domain (i.e. including non-trivial bundles) over a smooth reduced  curve is naturally identified with the Picard variety of bundles of a fixed degree. Here,  the spectral curve is an elliptic curve
 \cite{Hi}.

\subsection{Real structure}
 
 This of course neglects the real structure: the fact that we are considering the coadjoint orbits in ${\mathfrak{su}(3)}$ instead of ${\mathfrak{sl}(3,\CC)}$. This implies that $L(\overline z)^* = L(z)$; on the level of the spectrum, the spectral curve is invariant under 
 $$I(z,\eta) = (\overline z,-\overline \eta).$$ This invariance reduces the dimension of the set of spectral curves for our reduced coadjoint orbits from one complex dimension to one real. The real (invariant) portion of an elliptic curve has one or two components. Let us contrast this with the spectrum of $L(z)$ over the real line. The fact that elements of ${\mathfrak{su}(3)}$ are always diagonalisable tells us that multiple eigenvalues for a given $z$ correspond to a singularity of the spectral curve. If the curve is not to be singular, one has an ordering of ($i$ times) eigenvalues$$i\eta_1(z)<i \eta_2(z) <i\eta_3(z)$$over the real line. This would seem to give three real components to the spectral curve, except that the curve is sitting inside $\calo(1)$, whose real part is a Moebius strip, containing the real part of the curve. Compactifying joins two of the real components that we had over $\RR$ together so that while there are three components over  $\RR$, the real part of our compact spectral curve only has two components.

 On the level of line bundles, one needs a positive definite Hermitian form on $H^0(S,\calE)$; to obtain this, one must have $\calE\otimes \overline{I^*\calE} = K_S(2) = \calo(2)$ (\cite{AHH}). This reduces one to a real locus on the Jacobian, since any two $\calE, \calE'$ differ by a degree zero line bundle $L$ satisfying  $L\otimes \overline{I^*L}= \calo$; such bundles form a disjoint union of two circles in the Jacobian of $S$. These two real circles in the Jacobian    correspond to Hermitian forms of different signature, with one being definite, and the other indefinite; definiteness for the form imposes extra constraints, and reduces us to one circle for each spectral curve.
 
\section{Hamiltonian flow}

\subsection{Darboux coordinates}

 The general theory of Hitchin systems relates our Lax flows to a linear flow of line bundles along the curve, that is, along our real circles. There is an explicit method of integration, given in \cite{AHH2}, in terms of ``divisor coordinates''. Indeed, one can represent a line bundle on a curve by a divisor $\sum_i(z_i,\eta_i)$, and the main theorem of \cite{AHH2} is that this gives us Darboux coordinates, after a bit of normalisation.  In our case, there will just be one point $(z_0,\eta_0)$ in the  divisor.

For our case, the Darboux coordinates are fairly easy to compute.
 Indeed, one can conjugate so that $Y-Z$ is diagonal, with eigenvalues $\alpha_1, \alpha_2,\alpha_3$. Our Lax matrix is then 
 $$L(z) = \begin{pmatrix} \alpha_1&0&0\\ 0& \alpha_2&0\\ 0& 0&\alpha_3\end{pmatrix} z + \begin{pmatrix} s_{1,1}&s_{1,2}&s_{1,3}\\ s_{2,1}& s_{2,2}&s_{2,3}\\s_{3,1}& s_{3,2}&s_{3,3}\end{pmatrix}.$$
As in \cite{AHH2}, one can compute our coordinates by the vanishing of  a column (say the first) of the matrix $\widetilde L(z,\eta)$ of cofactors of $L(z)-\eta$; this follows from the relation $ (L(z)-\eta)\widetilde L(z,\eta) = \rho(z,\eta) I$ . The two non-diagonal elements of this are given by 
$$\widetilde L_{2,1} (z,\eta)= -s_{2,1}(\alpha_3 z-\eta + s_{3,3}) + s_{3,1} s_{2,3}$$
$$\widetilde L_{3,1} (z,\eta)= - s_{3,1}(\alpha_2 z-\eta + s_{2,2}) + s_{2,1} s_{3,2};$$
  that is, by two linear functions in $z,\eta$.  Setting both of them equal to zero, and solving, one obtains a unique solution
  $z_0,\eta_0$, and hence $z_0,\zeta_0$.  One then checks that $\widetilde L_{1,1} (z_0,\eta_0)$ also vanishes. Setting $\zeta_0 = \eta_0/z_0(z_0^2-1)$, the theorem is that $z_0,\zeta_0$ are Darboux coordinates.
  
\subsection{The flow}

  A classical Liouville generating function technique then gives the flow. The curve $S = S(H)$ defines implicitly $\eta $, or alternately, $\zeta = \eta/z(z^2-1) $, as a function of $z$, and $H$: setting
  $$G(z_0,H) =\int^{z_0} \zeta(z,H) dz$$
  so that $\frac{\partial G}{\partial z_0} = \zeta_0$, one has that the linearizing coordinate for the flow  is 
  $$t(z_0) = \frac {\partial G}{\partial H} = \int_0^{z_0}\frac{-\partial \rho/\partial H}{ z(z^2-1)\partial \rho/\partial \eta} dz= \int_0^{z_0} \frac{dz}{
Q_1(z)    - 3\eta^2}$$
with the integral taken along the spectral curve. The integrand, via 
the Poincar\'e residue theorem, is a holomorphic one-form  on the curve. As we are flowing on a circle, let us suppose that, starting out at $t(z_0) = 0$, one returns to the initial value at $t(z_0) = T$; a priori, this first return time depends on $H$: $T= T(H)$.  One then has for the sheaves $\calE(T(H)) = \calE(0)$, and so by Abel's theorem   the value of $T(H)$ is equal to the integral of the form around a cycle $\gamma$ on the curve, so that one has a period
$$T(H) = \int_\gamma \frac{dz}{
Q_1(z)    - 3\eta^2}.$$
Some function $F(H)$ should give the circle action, and this has vector field equal to  $\frac{dF}{dH}$ times that of $H$, and so a period that is $\frac{dF}{dH}^{-1}$ times that of $H$.  
One thus has the differential equation
$$T(H)   = \frac{dF}{dH}$$
for the Hamiltonian of the circle action. Returning to the definition of our linearizing variable,  a solution to this equation is given by 
$$F(H) = \int_\gamma \zeta(z,H) dz.$$
Of course, as $z, \zeta$ are Darboux coordinates, this is the technique given for example in the book of Arnol'd \cite{Arn} (Chap. 10, \S 50) for passing to action-angle variables. We would like to also note that the form $\zeta(z,H) dz$ plays a particular role in the Hitchin integrable systems, and more generally in integrable hierarchies --- see the paper of Krichever and Phong \cite{KP}.

\subsection {Properties of $H$, $F(H)$, and the curve}

Let us consider the possible range of $H$. One has the characteristic polynomial
 $\rho(z,\eta) = (iH z(z^2-1) + Q_0(z)) + Q_1(z) \eta   - \eta^3$, with imaginary roots $\eta_i$ along the real axis. Transforming the equation, so that the roots are real roots $i\eta_i$, one has a real polynomial 
$$ (H z(z^2-1)  -i Q_0(z)) - Q_1(z) \eta   - \eta^3.$$
Recall that the fact we are dealing with elements of ${\mathfrak{su}(3)}$ tells us that we have three real solutions to this equation for all $z$. 
The polynomial has discriminant$$-27 (H z(z^2-1)  -i Q_0(z))^2  - 4 Q_1(z)^3.$$ For the equation to have three distinct real roots for all $z$, one wants 
\begin{align*}  -27 (H z(z^2-1)&  -i Q_0(z))^2  - 4 Q_1(z)^3  = \\
&-27 (z(z^2-1))^2 H^2 -(54i z(z^2-1)Q_0(z))  H +27Q_0(z)^2 -4Q_1(z)^3 >0 \end{align*}
for all $z$. In particular, if $c$ is the leading ($z^6$) term of $-4 Q_1(z)^3$, one wants
$$ -27 H^2 + c >0,$$
showing that $H$ lies in an interval $I$; at the ends of this interval, the curve $S$ acquires a singularity at infinity.
 
 A further remark is that the graph  in $H$ of the  discriminant, for each fixed real $z\neq -1,0,1$,  is a downward-pointing parabola. Thus, if it is positive or zero for $H_1$ and $H_2$, it is positive for every $H$ in the interval $[H_1,H_2]$.  Hence, a good strategy is to take the curves for which $ -27 H^2 + c =0,$ and show that they have three real points over each finite real $z$. This would show that the spectral curves for $H$ with $ -27 H^2 + c >0$ have no singular points along the real line. 
 
 On the other hand, if they do correspond to   $L(z)$ with   values in  ${\mathfrak{su}(3)}$, they cannot have a singular point away from the real line: because of the real structure, such a singularity would come in pairs. A singularity corresponds to the vanishing of (all of) the matrix of cofactors, the off diagonal terms of which are   linear in the coordinates $\eta, z$ and so can have at most one solution. Another way of seeing that we can have only one singular point is to pass to the Weierstrass normal form of the curve, $y^2 = $ cubic in $x$, by a projective change of coordinates. Singular points of the curve then correspond to multiple zeroes of the cubic, and there can only be one.
 
 This then leads one to consider the curves for $ -27 H^2 + c =0$ which have a singularity at infinity. Again, this tells us that there cannot be a singularity elsewhere, and so it must be smooth at finite points. The matrices corresponding to $ -27 H^2 + c =0$, with two equal eigenvalues at infinity, can be normalized to 
 $$L(z) = \begin{pmatrix} -2\alpha &0&0\\ 0& \alpha &0\\ 0& 0&\alpha \end{pmatrix} z + \begin{pmatrix} s_{1,1}&s_{1,2}&s_{1,3}\\ s_{2,1}& s_{2,2}&0\\s_{3,1}& 0&s_{3,3}\end{pmatrix};$$
 one can check that for these the discriminant is a quartic, positive  in a neighbourhood of infinity. Since there are no other singular points, the discriminant remains positive over the whole line. For curves with  $ -27 H^2 + c >0$, then, one has a positive discrimant everywhere.
 
 Note that this picture of singular curves at the end of the interval fits well with the corresponding $L(z)$ being a fixed point of the Hamiltonian action; the reason is that the rank one sheaves on a singular nodal curve (the case we consider here) have a generic stratum corresponding to line bundles on the curve, and a codimension one stratum, corresponding to direct images of line bundles on the desingularization of the curve. In terms of the $L(z_0)$ at the corresponding point $z_0$ of the line, one has the generic $L(z_0)$, with a double eigenvalue but a generic non-diagonal Jordan form, and the codimension one stratum , with $L(z_0)$ diagonalizable. The fact that we are in ${\mathfrak{su}(3)}$ forces this non-generic case. In the case which concerns us, the singular curve has arithmetic genus one, and its blow-up is rational. Line bundles on the rational curves are rigid, and so we have a fixed point of the flow.

\subsection{Hamiltonian of $S^1$ action} 

 Finally, there remains the question of the Hamiltonian of the circle action. We have a differential equation for it, but one would like to be reassured that the Hamiltonian takes values in a compact interval, giving the standard picture of the Hamiltonian for a circle action as a standard coordinate function on the standard 3-sphere. This no doubt can be seen from the general theory of toric varieties, but can be seen directly as follows. We had the equations for our Darboux coordinates $z_0, \eta_0$, i.e.
 
 $$0= -s_{2,1}(\alpha_3 z_0-\eta_0 + s_{3,3}) + s_{3,1} s_{2,3}=  -s_{3,1}(\alpha_2 z_0-\eta_0 + s_{2,2}) + s_{2,1} s_{3,2},$$
 and so 
 $$ \alpha_3 z_0-\eta_0  = \frac {s_{3,1} s_{2,3}}{s_{2,1}} -  s_{3,3} $$
 $$\alpha_2 z_0-\eta_0 =  \frac{s_{2,1} s_{3,2}}{ s_{3,1}} - s_{2,2},$$
 giving 
 $$z_0 = \frac{1} { \alpha_2- \alpha_3} \left( \frac{s_{2,1} s_{3,2}}{ s_{3,1}} - s_{2,2}- \frac{s_{3,1} s_{2,3}}{s_{2,1}} + s_{3,3}\right).$$ This produces a coordinate that tends to infinity as  $ \alpha_2- \alpha_3$  tends to $0$ at $$-27 H^2 + c=0.$$ Accordingly, $T(H)$ is now an integral along a vanishing cycle on the curve, but of a form given as a Poincar\'e residue, which, at the limit, has a simple pole at the singular point with a non-zero residue. Thus, the reparametrization of $H$ stays finite at the end points, and the new Hamiltonian $F(H)$ takes values in an interval.\\

\noindent\emph{Remark.}  In the case of the $n$-fold reduced product in $\mbox{SU}(m)$, there is a corresponding picture within a moduli space of rank-$m$ Higgs bundles on the $n$-punctured sphere. One has spectral curves, and line bundles, and a real involution part of whose fixed point locus corresponds to  the reduced product. One also has an integrable system;  however, the passage to action-angle coordinates (i.e. extracting a torus action from the integrable system) is a lot less evident in higher dimensions, though the technique of passing to 
action-angle coordinates also works there. What is less certain is whether  this works well globally,  for example  when the action of the torus ceases to be free, e.g. at singular spectral curves.  Nevertheless, one can ask, after passing to a subset of real spectral curves, whether there is some sort of polytope's worth of spectral curves with, for example, singular curves at the boundary giving us something resembling a torus action. 

 We note from the paper \cite{JW} that having such a polytope does not always guarantee that one has a toric variety.

\end{document}